\documentstyle[twoside, amsmath, amssymb, amsfonts, eucal, epsfig,pstricks,pst-node,pst-tree,ifthen]{article}
\newtheorem{theorem}{Theorem}[section]

\newtheorem{proposition}[theorem]{Proposition}
\newtheorem{corollary}[theorem]{Corollary}
\newtheorem{definition}[theorem]{Definition}

\newtheorem{algorithm}[theorem]{Algorithm}

\newcount\refno
\newcommand\la{\langle}
\newcommand\ra{\rangle}

\newcommand\be{\begin{equation}}
\newcommand\ee{\end{equation}}
\newcommand\bn{\begin{eqnarray}}
\newcommand\en{\end{eqnarray}}
\newcommand\bns{\begin{eqnarray*}}
\newcommand\ens{\end{eqnarray*}}

\newcommand\bC{{\mathbb C}}
\newcommand\bR{{\mathbb R}}
\newcommand\bN{{\mathbb N}}

\newcommand\bZ{{\mathbb Z}}

\newcommand\bF{{\mathbb F}}

\newcommand\F{{\cal F}}

\newcommand\Ga{{\Gamma}}
\newcommand\ga{{\gamma}}
\newcommand\ep{{\epsilon}}

\newcommand{\NS}{\mbox{\scriptsize${\mathbb{N}}$}}

\def\matrice#1{\left(#1_{n,k}\right)_{n,k\in \NS}}

\def\divdif{\mathord\kern.43em\vrule width.6pt height5.6pt depth-.28pt \kern-.43em\Delta}

\title{
Generalized Stirling Numbers and Generalized Stirling Functions
} 

\author{ 
Tian-Xiao He\footnote{This paper was presented in the International Conference on Asymptotics and Special Functions, City University of Hong Kong, 
30 May - 03 June, 2011}\\
{\small Department of Mathematics and Computer Science }\\
{\small Illinois Wesleyan University}\\
{\small Bloomington, IL 61702-2900, USA}\\
}

\date{April, 2011}

\begin{document}

\maketitle
\setcounter{page}{1}
\pagestyle{myheadings}
\markboth{ 
T. X. He} 
{Generalized Stirling Functionss 
}

\begin{abstract}
\noindent
Here presented is a unified approach to Stirling numbers and their generalizations as well as generalized Stirling functions by using generalized factorial functions, $k$-Gamma functions, and generalized divided difference. Previous well-known extensions of Stirling numbers due to Riordan, Carlitz, Howard, Charalambides-Koutras, Gould-Hopper, Hsu-Shiue, Tsylova Todorov, Ahuja-Enneking,  and Stirling functions introduced by Butzer and Hauss, Butzer, Kilbas, and Trujilloet and others are included as particular cases of our generalization. Some basic properties related to our general pattern such as their recursive relations and generating functions  are discussed. Three algorithms for calculating the Stirling numbers based on our generalization are also given, which include a comprehensive algorithm using the characterization of Riordan arrays. 


\end{abstract}

\section{Introduction}
\setcounter{equation}{0}

The classical Stirling numbers of the first kind and the second kind, denoted by $s(n,k)$ and $S(n,k)$, respectively, can be defined via a pair of inverse relations 

\be\label{eq:1.1} 
[z]_n=\sum^n_{k=0} s(n,k) z^k,\quad z^n=\sum^n_{k=0} S(n,k) [z]_k,
\ee
with the convention $s(n,0)=S(n,0)=\delta_{n,0}$, the Kronecker symbol, where $ z\in {\bC}$, $n\in {\bN}_0={\bN}\cup\{0\}$, and the falling factorial polynomials $[z]_n=z(z-1)\cdots (z-n+1)$.  $|s(n,k)|$ presents the number of permutations of $n$ elements with $k$ disjoint cycles while $S(n,k)$ gives the number of ways to partition $n$ elements into $k$ nonempty subsets. The simplest way to compute $s(n,k)$ is finding the coefficients of the expansion of $[z]_n$. \cite{HS10} gives a simple way to evaluate $S(n,k)$ using Horner's method. 

Another way of introducing classical Stirling numbers is via their exponential generating functions

\be\label{eq:1.2}
\frac{(\log (1+x))^k}{k!} =\sum_{n\geq k} s(n,k) \frac{x^n}{n!}, \quad 
\frac{(e^x-1)^k}{k!} =\sum_{n\geq k} S(n,k) \frac{x^n}{n!},
\ee
where $|x|<1$ and $k\in {\bN}_0$. In \cite{Jor}, Jordan said that, ``Stirling's numbers are of the greatest utility. This however has not been fully recognized.'' He also thinks that, ``Stirling's numbers are as important or even more so than Bernoulli's numbers.'' 

Besides the above two expressions, the Stirling numbers of the second kind has the following third definition (see \cite{Com} and \cite{Jor}), which is equivalent to the above two definitions but makes a more important rule in computation and generalization. 

\bn\label{eq:1.3}
S(n,k)&:=&\frac{1}{k!} \left. \Delta^k z^n\right|_{z=0}=\frac{1}{k!} \sum^k_{j=1} (-1)^{k-j} {k\choose j} j^n\nonumber\\
&=&\frac{1}{k!}\sum^k_{j=1} (-1)^j{k\choose j} (k-j)^n.
\en
Expressions (\ref{eq:1.1}) - (\ref{eq:1.3}) will be our starting points to extend the classical Stirling number pair and the Stirling numbers.

Denote $\la z\ra_{n,\alpha}:= z(z+\alpha)\cdots (z+(n-1)\alpha)$ for $n=1,2,\ldots$, and $\la z\ra_{0,\alpha}=1$, where $\la z\ra_{n,\alpha}$ is called the generalized factorial of $z$ with increment $\alpha$. Thus, $\la z\ra_{n,-1}=[z]_n$ is the classical falling factorial with $[z]_0=1$, and $\la z\ra_{n,0}=z^n$. More properties of $\la z\ra_{n,\alpha}$ will be presented below. 

With a closed observation, Stirling numbers of two kinds defined in (\ref{eq:1.1}) can be written as a unified Newton form:

\be\label{eq:1.1*}
\la z\ra_{n,-\alpha} =\sum^n_{k=0} S(n,k,\alpha, \beta)\la z\ra_{n,-\beta}, 
\ee
with $S(n,k,1,0)=s(n,k)$, the Stirling numbers of the first kind and $S(n,k,0,1)=S(n,k)$. the Stirling numbers of the second kind. 
Inspired by (\ref{eq:1.1*}) and many extensions of classical Stirling numbers or Stirling number pairs introduced by \cite{Car80}, \cite{Hor85}, \cite{Tsy}, \cite{HS98}, etc. We may define a unified {\it generalized Stirling numbers} $S(n,k,\alpha, \beta,r)$ as follows. 

\begin{definition}\label{def:1.0} 
Let $n\in{\bN}$ and $ \alpha, \beta, r\in{\bR}$. A generalized Stirling number 
denoted by $S(n,k,\alpha, \beta, r)$ is defined by 

\be\label{eq:1.1**}
\la z\ra_{n,-\alpha}=\sum^n_{k=0} S(n,k,\alpha,\beta,r)\la z-r\ra_{k,-\beta}.
\ee
In particular, if $(\alpha, \beta,r)=(1,0,0)$, $S(n,k,1,0,0)$ is reduced to the unified form of Classical Stirling numbers defined by (\ref{eq:1.1*}).
\end{definition}

Each $\la z\ra_{n,-\alpha}$ does have exactly one such expansion (\ref{eq:1.1**}) for any given $z$. Since $deg\,\, \la z-r\ra_{k,-\beta}=k$ for all $k$, which generates a graded basis for $\Pi\subset {\bF}\to {\bF}$, the linear spaces of polynomials in one real (when ${\bF}={\bR}$) or complex (when ${\bF}={\bC}$), in the sense that, for each $n$, $\{ \la z-r\ra_{n, -\beta}\}$ is a basis for $\Pi_n\subset \Pi$, the subspace of all polynomials of degree $<n$. In other wards, the column map 

\[
W_z: {\bF}^N_0 \to \Pi: s\mapsto \sum_{k\geq 0} S(n,k,\alpha, \beta, r) \la z\ra_{k,-\beta},
\]
from the space ${\bF}^N_0$ of scalar sequences with finitely many nonzero entries to the space $\Pi$ is one-to-one and onto, hence invertible. In particular, for each $n\in {\bN}$, the coefficient $c(n)$ in the Newton form (\ref{eq:1.1**}) for $\la z\ra_{n,-\alpha}$ depends linearly on $\la z\ra_{n,-\alpha}$, i.e., $\la z\ra_{n,-\alpha}\mapsto s(n)=(W^{-1}_z \la z\ra_{n,-\alpha}) (n)$, the set of $S(n,k,\alpha, \beta,r)$, is a well-defined linear functional on $\Pi$, and vanishes on $\Pi_{< n-1}$. 

Similarly to (\ref{eq:1.1}), from Definition \ref{def:1.0} a Stirling-type pair $\{ S^1, S^2\} =\{ S^1(n,k),$ $ S^2(n,k)\}\equiv \{ S(n,k;$ $\alpha, \beta, r),$ $S(n,k; \beta, \alpha, -r)\}$ (see also in \cite{HS98}) can be defined by the inverse relations

\bn\label{eq:1.4}
\la z\ra_{n,-\alpha}&=&\sum^n_{k=0} S^1(n,k) \la z-r\ra_{k,-\beta}\nonumber\\
\la z\ra_{n,-\beta}&=& \sum^n_{k=0} S^2(n,k) \la z+r\ra_{k,-\alpha},
\en
where $n\in {\bN}$ and the parameter triple $(\alpha, \beta, r)\not= (0,0,0)$ is in ${\bR}^3$ or ${\bC}^3$. Hence, we may call $S^1$ and $S^2$ an $(\alpha, \beta,r)$
and a $(\beta, \alpha, -r)-$ pair. Obviously, 

\[
S(n,k; 0,0,1)={n\choose k}
\]
because $z^n=\sum^n_{k=0} {n\choose k} (z-1)^k$. In addition, the classical Stirling number pair $\{ s(n,k), S(n,k)\}$ is the $(1,0,0)-$ pair $\{ S^1, S^2\}$, namely,

\[
s(n,k)=S^1(n,k; 1,0,0)\quad S(n,k)=S^2(n,k; 1,0,0).
\]
For brevity, we will use $S(n,k)$ to denote $S(n,k,\alpha, \beta,r)$ if there is no need to indicate $\alpha$, $\beta$, and $r$ explicitly. From (\ref{eq:1.1**}), one may find 

\be\label{eq:1.5}
S(0,0)=1,\quad S(n,n)=1,\quad S(1,0)=r, \quad and\quad S(n,0)=\la r\ra _{n,-\alpha}.
\ee
Evidently, substituting $n=k=0$ into 
(\ref{eq:1.1**}) yields the first formula of (\ref{eq:1.5}). Comparing the coefficients of the highest power terms on the both sides of (\ref{eq:1.1**}), we obtain the second formula of (\ref{eq:1.5}).  Let $n=1$ in (\ref{eq:1.1**}) and noting $S(1,1)=1$, we have the third formula. Finally, substituting $z=r$ in (\ref{eq:1.1**}), one can establish the last formula of (\ref{eq:1.5}). The numbers $\sigma (n,k)$ discussed by Doubilet et al. in \cite{DRS} and by Wagner in \cite{Wag} is $k! S(n,k;0,1,0)$. More special cases of the generalized Stirling numbers and Stirling-type pairs defined by (\ref{eq:1.1**}) or (\ref{eq:1.4}) are surveyed below in Table 1. 

\begin{center}
\begin{tabular}{c|c|c|c}
$  (\alpha, \beta, r)$ & $ dual\,\, of\,\, (\alpha, \beta, r) $ 
&$  S(n,k) \,\,pairs $ &$  Name \;\; of \;\; Stirling\,\, numbers $\\ \hline \rule[-3mm]{0mm}{8mm}
$(-1,1,0) $  & $(1,-1,0)$ & $ \begin{array}{cc} n!{n-1\choose {k-1}}/k! \\
(-1)^{n-k}n! {n-1\choose{k-1}}/k! \end{array}$  & $Lah\ number \ pair \cite{JRS}$\\ \hline
\rule[-3mm]{0mm}{8mm}
$(-1,0,0)$  & $(0,-1,0)$ 
&   $\begin{array}{cc}|s(n,k)|\\ (-1)^{n-k} S(n,k)\end{array}$  & $signless \ Stirling \ numbers\cite{Rot}$\\ \hline
\rule[-3mm]{0mm}{8mm}
$(1,\theta,0) (\theta\not= 0)$ & $(\theta,1,0)$ &  $\begin{array}{cc}S(n,k,1,\theta,0) \\ S(n,k, \theta, 1,0)
\end{array}$ & $\begin{array}{l}Carlitz's\ degenerate\ Stirling \\
number \ pair\cite{Car79}\end{array}$\\ \hline 
\rule[-3mm]{0mm}{8mm}
$(1,0,-\lambda)$  & $(0,1,\lambda)$
& $\begin{array}{cc} S(n,k,1,0,-\lambda)\\ S(n,k, 0,1,\lambda)\end{array}$ & $\begin{array}{l}Carlitz's \ weighted\ Stirling\\
number \ pair \cite{Car80}\end{array}$\\ \hline
\rule[-3mm]{0mm}{8mm}
$(1,\theta, -\lambda)$      & $(\theta, 1, \lambda)$ &  $\begin{array} {cc} S(n,k,1,\theta, -\lambda)\\S(n,k, \theta, 1, \lambda)\end{array}$ 
& $\begin{array} {l}Howard's\ weighted\ degenerate\\ Stirling 
\ number\  pair \cite{Hor85}\end{array}$\\ \hline
\rule[-3mm]{0mm}{8mm}
$(0,1,-a+b)$  & $(1,0, -b+a)$        &  $\begin{array}{cc}S(n,k,0,1,-a+b)\\S(n,k,1,0,-b+a)\end{array}$& 
$\begin{array}{l} Gould-Hopper's\ non-central\ Lah\\ number \ pair \cite{GH}\end{array}$
\footnotemark\\ \hline
\rule[-3mm]{0mm}{8mm}
$(1/s,1,-a+b)$ & $(1,1/s, -b+a)$ & $\begin{array}{cc}S(n,k,1/s,1,-a+b)\\S(n,k, 1,1/s, -b+a)\end{array}$ & 
$\begin{array} {l} Charalambides-Koutras's\ non-\\central \ C \ number \ pair\cite{Cha,CK}\end{array}$\footnotemark 
\\ \hline
\rule[-3mm]{0mm}{8mm}
$(1,0, b-a)$  & $(0,1,a-b)$ &  $\begin{array}{cc}S(n,k,1,0,b-a)\\S(n,k,0,1,a-b)\end{array}$ & 
$\begin{array}{l} Riordan's\ non-central \ Stirling\\
number\ pair \cite{Rio}\end{array}$ \\ \hline
\rule[-3mm]{0mm}{8mm}
$(\alpha,\beta,0)$ & $(\beta,\alpha,0)$ & 
$\begin{array}{cc}A_{\alpha\beta}(r,m)\\B_{\alpha\beta}(r,m)\end{array}$& $Tsylova's\ Stirling\ number\ pair\cite{Tsy}$\\ \hline
\rule[-3mm]{0mm}{8mm}
$(\alpha,\beta,r)$ & $(\beta,\alpha,-r)$
& $\begin{array}{cc}S(n,k, \alpha, \beta,r)\\S(n,k,\beta, \alpha,-r)\end{array}$ & $\begin{array}{l} Hsu-Shiue's \ Stirling\\ number\ pair\cite{HS98}\end{array}$\\ \hline
\rule[-3mm]{0mm}{8mm}
$(1,x,0)$ & 
-- &       
$a_{nk}(x)$  & $Todorov's\ Stirling\ numbers\cite{Tod}$  \\ \hline
\rule[-3mm]{0mm}{8mm}
$(-1/r,1,0)$ & 
 -- &  $B(n,r,k)$ & $\begin{array}{l} Ahuja-Enneking's
\ associated\\  Lah\ numbers\cite{ND}\footnotemark \end{array}$
\\ \hline
\rule[-3mm]{0mm}{8mm}
$(-1,0,r)$  & 
-- &   $S(n-r, k-r, -1,0,r)$ & $Broder's \ r-Stirling\ numbers \cite{Bro}$\\ \hline
\rule[-3mm]{0mm}{8mm}
\end{tabular}
\centerline{Table 1. Some generalized Stirling Numbers and Stirling Number pairs}
\end{center}
\vspace{.15in}

The classical falling factorial polynomials $[z]_n=z(z-1)\cdots (z-n+1)$ and classical rising factorial polynomials $[z]^n=z(z+1)\cdots (z+n-1)$, $z\in{\bC}$ and $n\in{\bN}$, can be unified to the expression 

\[
\la z\ra_{n,\pm 1}:=z(z\pm 1)\cdots (z\pm (n-1)),
\]
using the {\it generalized factorial polynomial} expression

\be\label{eq:1.6}
\la z\ra_{n,k}:=z(z+k)\cdots (z+(n-1)k)=\la z+(n-1)k\ra_{n,-k}\quad (z\in {\bC}, n\in {\bN}).
\ee
Thus $\la z\ra_{n,1}=[z]^n$ and $\la z\ra_{n,-1}=[z]_n$. In addition, we immediately have the relationship between $[z]^n$ and $\la z \ra_{n,k}$ as 

\be\label{eq:1.7}
\la z\ra_{n,k}=k^n[z/k]^n \quad (z\in {\bC}, n\in {\bN}, k> 0).
\ee
Similarly, we obtain 

\be\label{eq:1.8}
\la z\ra_{n,-k}=z(z-k)\cdots (z-(n-1)k)=k^n[z/k]_n\quad (z\in {\bC}, n\in {\bN}, k> 0).
\ee 
The history as well as some important basic results of the generalized factorials can be found in Chapter II of \cite{Jor}, and an application of the generalized factorials in the Lagrange interpolation is shown on Page 31 of \cite{Gel}.

It is known that the falling factorial polynomials and rising factorial polynomials can be presented in terms of Gamma functions: $[z]_n=\Ga(z+1)/\Ga(z-n+1)$ and $[z]^n=\Ga (z+n)/\Ga(z)$, and the gamma function $\Gamma(z)$ can be defined in terms of factorial functions by (see, for example, \cite{MK})

\be\label{eq:1.9} 
\Ga(z)=\displaystyle \lim_{n\to \infty} \frac{n! n^{z-1}}{[z]^n} \quad (z\in {\bC}-k{\bZ}_-). 
\ee
As an analogy, the {\it $k$-gamma function} $\Ga_k$, a one parameter deformation of the classical gamma function, is defined by (see, for example \cite{DP}) 

\be\label{eq:1.10}
\Ga_k (z):=\displaystyle\lim_{n\to \infty} \frac{n! k^n (nk)^{\frac{z}{k}-1}}{\la z\ra_{n,k}}\quad (k> 0, z\in {\bC}-k{\bZ}_-).
\ee
$[z]^n$ and $\la z\ra_{n,k}$ ($k>0$) are also called the {\it Pochhammer symbol} and {\it $k$-Pochhammer symbol,} respectively. Even the parameter $k$ is replaced by 
other parameters, we still call the corresponding Pochhammer symbol the $k$-Pochhammer. 

For $k>0$, from (\ref{eq:1.7}), (\ref{eq:1.9}) and (\ref{eq:1.10}) (see also \cite{Man}) we have

\be\label{eq:1.11}
\Ga_k (z)=k^{(z/k)-1} \Ga\left( \frac{z}{k}\right).
\ee

Since $[z]^n=\Ga(z+n)/\Ga (z )$, \cite{BH93} extends the classical raising and falling factorial polynomials to generalized raising and falling functions associated with real number $\ga$ by setting 

\be\label{eq:1.12}
[z]^\ga:=\frac{\Ga(z+\ga)}{\Ga (z)}\quad [z]_\ga :=\frac{\Ga(z+1)}{\Ga(z-\ga +1)},
\ee
respectively. We now extend $\la z\ra_{n,k}$ defined by (\ref{eq:1.6}) to a generalized form associated with $\ga\in{\bC}$ using the relationship (\ref{eq:1.7}), namely,

\be\label{eq:1.13}
\la z\ra_{\ga,k}=k^\ga [z/k]^\ga,\quad  \la z\ra_{\ga, -k} =k^\ga [z/k]_\ga \quad (z\in {\bC}, \ga \in {\bC}, k>0),
\ee
which are called the generalized raising and falling factorial functions associated with complex number $\ga$, respectively. Using (\ref{eq:1.11})-(\ref{eq:1.13}), we establish the following result. 

\begin{theorem}\label{thm:1.1}
If $k>0$ and $\la z\ra_{\ga ,k}$ is defined by (\ref{eq:1.13}), then 

\be\label{eq:1.14}
\la z\ra_{\ga ,k}=\frac{\Ga_k(z+\ga k)}{\Ga_k (z)}\qquad \la z\ra_{\ga, -k} =\frac{\Ga_k(z+k)}{\Ga_k(z-(\ga-1)k)}.
\ee
\end{theorem}

There hold the following recurrence relations of  the generalized raising and falling factorial functions. 

\begin{proposition}\label{pro:2.2}
If $k>0$ and $\la z\ra_{\ga ,k}$ is defined by (\ref{eq:1.13}), then there hold 

\be\label{eq:1.17}
\la z\ra_{\ga ,k}=(z+(\ga -1)k)\la z\ra_{\ga-1,k}, \quad \la z\ra_{\ga, -k} =(z-(\ga-1)k)\la z\ra_{\ga-1,-k}.
\ee
\end{proposition}

In next section, we will present the unified expression and some properties of the generalized Stirling numbers of integer orders. Two algorithms based on the unified expression will be given. Then, we use the $k$-Pochhammer symbol and $k$-Gamma functions to extend the classical Stirling numbers of integer orders to the complex number orders in Section $3$, which are called the {\it generalized Stirling functions}. The convergence and the recurrence relation of the generalized Stirling functions as well as their generating functions will also be presented. Finally, in Section $4$ we will give more properties of generalized Stirling numbers and functions using the generating functions of generalized Stirling numbers shown in Section $3$, which include the asymptotic expansions of generalized Stirling numbers and functions and the sequence characterizations of the Riordan arrays of generalized Stirling numbers. The third algorithm of the computation of the generalized Stirling numbers, including the classical Stirling numbers as a special case, will be shown using the characterizations of their Riordan arrays.    

\section{Expressions of generalized Stirling numbers}
\setcounter{equation}{0}

First, we give an equivalent form of the generalized Stirling numbers $S(n,k)$ defined by (\ref{eq:1.1**}) by using the {\it generalized difference operator} 
in terms of $\beta$ ($\beta\not= 0$) defined by 

\be\label{eq:1.15*}
\Delta^k_\beta f =\Delta_\beta(\Delta^{k-1}_\beta f)\quad  (k\geq 2) \quad and \quad \Delta_\beta f(t):=f(t+\beta)-f(t).
\ee 
It can be seen that $\left. \Delta^k_\beta \la z\ra_{j,-\beta}\right|_{z=0}=\beta^k k!\delta_{k,j}$, where $\delta_{k,j}$ is the Kronecker delta symbol; i.e., $\delta_{k,j}=1$ when $k=j$ and $0$ otherwise. Evidently, from (\ref{eq:1.8}) there holds 

\be\label{eq:1.15-2}
\left. \Delta_\beta^k \la z\ra_{j, -\beta}\right|_{z=0}=\left. \Delta^k_\beta \beta^j\left[ \frac{t}{\beta}\right]_j\right|_{z=0}=\left. \beta^j\Delta^k [t]_j\right|_{z=0}=\beta^k k! \delta_{k,j}.
\ee

Denote the divided difference of $f(t)$ at $t+i$, $i=0,1,\ldots, k$, by $f[t, t+1, \ldots, t+k]$, or $[t,t+1, \ldots, t+k] f(t)$. Using the well-known forward difference formula, it is easy to check that 

\[
\frac{1}{k!}\Delta^k f(t)=f[t,t+1,\ldots, t+k]=[t,t+1,\ldots, t+k] f(t)
\]
and 

\[ 
\frac{1}{\beta^k k!}\Delta^k_\beta f(t)=f[t,t+\beta, t+2\beta,\ldots, t+k \beta]=[t,t+\beta, \ldots, t+k\beta] f(t).
\]
We now give the following definition of the {\it generalized divided differences.}

\begin{definition}\label{def:2.0}
We define ${\underline\triangle}^k_\beta f(t)$ by 

\be\label{eq:1.15**}
\underline \triangle^k_\beta f(t) =\left\{ \begin{array}{ll} \frac{1}{\beta^k k!}\Delta^k_\beta f(t)= f[t, t+\beta, \ldots, t+k \beta] 
& if\,\, \beta\not= 0\\ 
\frac{1}{k!}D^k f(t) &if\,\, \beta =0\end{array},\right.
\ee
where $\Delta^k_\beta f(t)$ is shown in (\ref{eq:1.15*}), $f[t,t+\beta, \ldots, t+k\beta]\equiv [t, t+\beta,\ldots, t+k\beta]f$ is the $k$th divided difference of $f$ in terms of $\{ t, t+\beta, \ldots, t+k\beta\}$, and $D^k f(t)$ is the $k$th derivative of $f(t)$.  
\end{definition}

From the well-known formula

\[
f[t,t+\beta, t+2\beta,\ldots, t+k \beta]=\frac{D^k f(\xi)}{k!},
\]
where $\xi$ is between $t$ and $t+k\beta$, it is clear that 

\be\label{eq:1.15***}
D^k f(t)=\lim_{\beta\to 0} \frac{1}{\beta^k} \Delta^k_\beta f(t),
\ee
which shows the generalized divided difference is well defined. 

We now give a unified expression of the generalized Stirling numbers in terms of the the generalized divided differences. 

\begin{theorem}\label{thm:2.1}
Let $n,k\in{\bN}_0$ and the parameter triple $(\alpha, \beta, r)\not= (0,0,0)$ is in ${\bR}^3$ or ${\bC}^3$. For the generalized Stirling numbers defined by (\ref{eq:1.1**}), there holds

\bn\label{eq:1.15-4}
&&S(n,k,\alpha,\beta,r)=\left. \underline \triangle^k_\beta \la z\ra_{n,-\alpha}\right|_{z=r} \nonumber\\
&=&\left\{ \begin{array}{ll} \frac{1}{\beta^k k!}\left. \Delta^k_\beta \la z\ra_{n,-\alpha}\right|_{z=r} 
=[r,r+\beta, \ldots, r+k\beta]\la z\ra_{n,-\alpha}& if\,\, \beta\not= 0\\
\frac{1}{k!}\left.  D^k \la z\ra_{n,-\alpha}\right|_{z=r} &if \,\, \beta=0. \end{array}\right.
\en
In particular, for the generalized Stirling number pair defined by (\ref{eq:1.4}), we have the expressions 

\bn\label{eq:1.15}
&&S^1(n,k)\equiv S^1(n,k,\alpha, \beta, r)=\left. \underline \triangle^k_\beta \la z\ra_{n,-\alpha}\right|_{z=r}
\nonumber\\
&=&\left\{ \begin{array} {ll} \frac{1}{\beta^k k!} \left. \Delta^k_\beta \la z\ra_{n,-\alpha}\right|_{z=r}=[r,r+\beta, \ldots, r+k\beta]\la z\ra_{n,-\alpha}, &if\,\,\beta\not= 0\\ 
\left. \frac{1}{k!}D^k\la z\ra_{n,-\alpha}\right|_{z=r}, &if\,\,\beta =0\end{array}\right.\\
&&S^2(n,k)\equiv S^2(n,k,\beta, \alpha,-r)=\left. \underline\triangle^k_\alpha \la z\ra_{n,-\beta}\right|_{z=-r}\nonumber\\
&=&\left\{ \begin{array}{ll}\frac{1}{\alpha^k k!}\left. \Delta^k_\alpha \la z\ra_{n,-\beta}\right|_{z=-r}
=[-r,-r+\alpha,\ldots, -r+k\alpha]\la z\ra_{n,-\beta}, & if\,\, \alpha\not= 0\\ \left. \frac{1}{k!}D^k\la z\ra_{n,-\beta}\right|_{z=-r}, &if\,\, \alpha =0\end{array}\right.\nonumber\\
\en
Furthermore, if $(\alpha, \beta,r)=(1,0,0)$, then (\ref{eq:1.15-4}) is reduced to the classical Stirling numbers of the first kind defined by (\ref{eq:1.1}) with the expression 

\[
s(n,k)=S(n,k,1,0,0)=\frac{1}{k!} \left. D^k [z]_n\right|_{z=0}.
\]
If $(\alpha,\beta,r)=(0,1,0)$, then (\ref{eq:1.15-4}) is reduced to the classical Stirling numbers of the second kind shown in (\ref{eq:1.3}) with the following divided difference expression form: 

\be\label{eq:1.16*}
S(n,k)=S(n,k, 0,1,0)=\left. [0,1,2,\ldots, k] z^n\right|_{z=0}.
\ee 
\end{theorem}

The following corollary is obvious due to the expansion formula of the divided differences generated from their definition.

\begin{corollary}\label{cor:2.2}
Let $n,k\in{\bN}_0$ and the parameter triple $(\alpha, \beta, r)\not= (0,0,0)$ is in ${\bR}^3$ or ${\bC}^3$. If $\beta\not= 0$, for the generalized Stirling numbers defined by (\ref{eq:1.1**}), there holds

\be\label{eq:1.16**}
S(n,k)\equiv S(n,k,\alpha,\beta,r )=\frac{1}{\beta^k k!}\sum^n_{j= 0} (-1)^j {k\choose j}\langle r+(k-j)\beta \rangle_{n,-\alpha}\quad (n\not= 0),
\ee
and $S(0,k)=\delta_{0k}.$
\end{corollary}

\noindent{\bf Remark 2.1} It can be seen from (\ref{eq:1.16**}) that 

\be\label{eq:1.16***}
S(n,0)\equiv S(n,0,\alpha,\beta, r)=\la r\ra_{n,-\alpha},
\ee
which is independent of $\beta$ and has been shown in (\ref{eq:1.5}). The difference is deriving (\ref{eq:1.16***}) from (\ref{eq:1.16**}) needs $(\alpha,r)\not= (0,0)$ when $\beta=0$. However, we have seen from (\ref{eq:1.5}) that the condition is not necessary. Another way to derive (\ref{eq:1.16***}) using the characterization of the Riordan arrays of the generalized Stirling numbers will be presented in the Algorithm \ref{alg:4.1} in Section $4$. 

\noindent{\bf Remark 2.2} If $\alpha\beta\not= 0$, by taking the $n$th forward differences in terms of $\alpha$ and $\beta$ on the both sides of two equations of (\ref{eq:1.4}), respectively, one may obtain identities 

\bns
n! \alpha^n &=&\sum^n_{k=0} S^1(n,k) \left. \Delta^n_\alpha \la z-r\ra_{k,-\beta}\right|_{z=0}\\
n! \beta^n &=& \sum^n_{k=0} S^2(n,k) \left. \Delta^n_{\beta} \la z+r\ra_{k,-\alpha}\right|_{z=0}.
\ens 
The above two identities can be unified to be one:

\[
n! \alpha^n =\sum^n_{k=0} \left. S(n,k,\alpha, \beta, r) \Delta^n_{\alpha} \la z-r\ra_{k,-\beta}\right|_{z=0}.
\]
When $\alpha=0$, the above identity turns to 

\[
n! =\sum^n_{k=0} S(n,k,0,\beta,r) \left. D^n \la z-r\ra_{k,-\beta}\right|_{z=0}.
\]

\noindent{\bf Remarker 2.3} 
There exists another expression of the divided difference $\left.\underline\triangle^k_\beta \la z\ra_{n,-\alpha}\right|_{z=r}$ in terms of {\it Peano kernel of B-spline}. Assume that the set $\tau:=\{ t, t+\beta,\ldots, t+k\beta\}$ lies in the interval $[a,b]$. Then on the interval, we have Taylor's identity

\[
\la z\ra_{n,-\alpha}=\sum_{j<k} \frac{(z-a)^j}{j!}\left. D^j \la z\ra_{n,-\alpha}\right|_{z=a}+\int^b_a\frac{(x-y)^{k-1}_+}{(k-1)!}  
\la y\ra_{n,-\alpha} dy.
\] 
If $\beta >0$, then $\underline\triangle^k_\beta$ is a weighted sum of values of derivatives of order $<k$, hence commutes with the integral in the above Taylor's expansion, which annihilates any polynomial of degree $< k$. Therefore, 

\[
\left. \underline\triangle^k_\beta \la z\ra_{n,-\alpha}\right|_{z=r} =\int^b_a \frac{M(y| \tau)}{k!} \la y\ra_{n,-\alpha} dy,
\]
where 

\[
M(y|\tau):=k [r,r+\beta,\ldots, r+k\beta] (\cdot -y)^{k-1}_+
\]
is the Curry-Schoenberg B-spline (see \cite{CS}) with the knot set $\tau$ and normalized to have integral $1$. In particular, 

\[
S(n,n,\alpha, \beta,r)=\left. \underline\triangle^k_\beta \la z\ra_{n,-\alpha}\right|_{z=r}  =\int^b_aM(y| r, r+\beta, \ldots, r+n\beta) dy =1.
\]

\medbreak
We now present two algorithms for calculating generalized Stirling numbers. If $\beta\not=0$, we denote 

\be\label{eq:no1}
\underline\triangle^j_\beta f(t+\ell\beta):=f[t, t+\ell \beta, t+(\ell +1)\beta, \ldots, t+j\beta]
\ee
Thus, from (\ref{eq:1.15-4}) in Theorem \ref{thm:2.1}, based on the recursive definition of the divided difference 
with respect to $\beta$ (see Definition \ref{def:2.0})

\be\label{eq:no2}
\underline\triangle^j_\beta f(t+\ell \beta)=\frac{1}{j\beta}(\underline\triangle^{j-1}_\beta f(t+(\ell +1)\beta)-\underline\triangle^{j-1}_\beta f(t+\ell\beta)),
\ee
we obtain an algorithm shown below. 

\begin{algorithm}\label{alg:1.1}
This algorithm of evaluating the generalized Stirling numbers is based on the construction of the following lower triangle array by using (\ref{eq:no1}) and (\ref{eq:no2}).

\begin{center}
\[
\begin{array}{lllll}
\left. 
\la z\ra_{n,-\alpha}\right|_{z=r} & & & & \\
\left. \la z+\beta \ra_{n,-\alpha} \right|_{z=r}& 
\left. \underline\triangle_{\beta}\la z\ra_{n,-\alpha} \right|_{z=r}& & &\\
\left. \la z+2\beta \ra_{n,-\alpha} \right|_{z=r}& 
\left. \underline\triangle_{\beta}\la z+\beta \ra_{n,-\alpha} \right|_{z=r}& 
\left. \underline\triangle^2_{\beta}\la z\ra_{n,-\alpha}\right|_{z=r} & & \\
\vdots &\vdots &\vdots &\ddots &\\
\left. \la z+k\beta \ra_{n,-\alpha}\right|_{z=r} &
\left. \underline\triangle_{\beta} \la z+(k-1)\beta \ra_{n,-\alpha}\right|_{z=r} &
\left. \underline\triangle^2_{\beta} \la z+(k-2)\beta\ra_{n,-\alpha}\right|_{z=r} & 
\cdots & \left. \underline\triangle^k_\beta \la z\ra_{n,-\alpha}\right|_{z=r} 
\end{array}
\]
\centerline{Table 2. 
The generalized Stirling numbers 
}
\end{center}
\vspace{.15in}
Thus, the diagonal of the above lower triangle array gives $S(n,i,\alpha,\beta,r)=\left.\underline\triangle^i_\beta\la z\ra_{n,-\alpha}\right|_{z=r}$ for $i=0,1,\ldots, k$.
\end{algorithm}

\noindent{\bf Example 2.1} We now use Algorithm \ref{alg:1.1} shown in Table $2$ to evaluate the classical Stirling numbers of the second kind $S(4,k)=S(4, k, 0,1,0)$ ($k=1,2,3,4$), which are re-expressed by (\ref{eq:1.16*}). Thus, 

\[
\begin{array}{lllll}
0 & \quad &\quad  & \quad & \quad  \\
1 &1 \quad & & & \\
2^4=16 & 15 \quad & 7 \quad & & \\
3^4=81 & 65 \quad & 25\quad  & 6\quad  &\\
4^4=256 & 175 \quad & 55 \quad & 10 \quad &1\quad
\end{array}
\]
From the diagonal of the above lower triangular matrix, we may read $S(4,0)=0$, $S(4,1)=1$, $S(4,2)=7$, $S(4,3)=6$, and $S(4,4)=1$. Meanwhile, the subdiagonal gives $S(5,1)=1$, $S(5,2)=15$, $S(5,3)=25$, and $S(5,4)=10$. 
 
 \medbreak

\noindent{\bf Example 2.2} For the Howard's weighted degenerate Stirling numbers $S(4,k)=S(4,k,1,1,-1)$, a similar argument of Example 2.1 yields 

\[
\begin{array}{lllll}
\left. \la z\ra_{4,-1}\right|_{z=-1}=24 & \quad &\quad  & \quad & \quad  \\
\left. \la z+1\ra_{4,-1}\right|_{z=-1}=0 &-24 \quad & & & \\
\left. \la z+2\ra_{4,-1}\right|_{z=-1}=0 & 0 \quad & 12 \quad & & \\
\left. \la z+3\ra_{4,-1}\right|_{z=-1}=0 & 0 \quad & 0\quad  & -4\quad  &\\
\left. \la z+4\ra_{4,-1}\right|_{z=-1}=0 & 0 \quad & 0 \quad & 0 \quad &1\quad
\end{array}
\]
Thus, $S(4,0)=24$, $S(4,1)=-24$, $S(4,2)=12$, $S(4,3)=-4$, and $S(4,4)=1$. 

\medbreak

\noindent{\bf Example 2.3} For the Howard's weighted degenerate Stirling numbers $S(4,k)=S(4,k,1,2,-1)$, using Algorithm \ref{alg:1.1}, we obtain $S(4,0)=24$, $S(4,1)=-12$, $S(4,2)=3$, $S(4,3)=2$, and $S(4,4)=1$ reading from the following table. 

\[
\begin{array}{lllll}
\left. \la z\ra_{4,-1}\right|_{z=-1}=24 & \quad &\quad  & \quad & \quad  \\
\left. \la z+2\ra_{4,-1}\right|_{z=-1}=0 &-12 \quad & & & \\
\left. \la z+4\ra_{4,-1}\right|_{z=-1}=0 & 0 \quad & 3 \quad & & \\
\left. \la z+6\ra_{4,-1}\right|_{z=-1}=120 & 60 \quad & 15\quad  & 2\quad  &\\
\left. \la z+8\ra_{4,-1}\right|_{z=-1}=840 & 360 \quad & 75 \quad & 10 \quad &1\quad
\end{array}
\]
\medbreak

\noindent{\bf Remark 2.4} Obviously, Algorithm \ref{alg:1.1} is not limited to the case of $\beta\not= 0$ since when $\beta =0$, $\left. \underline\triangle^k_{\beta}\la z\ra_{n,-\alpha}\right|_{z=r}$ ($k=0,1,\ldots, n$) on the diagonal of the lower triangle matrix in Table 1 are simply the $1/k!$ multiply of the derivatives $\left. D^k \la z\ra_{n,-\alpha}\right|_{z=r}$ (see Theorem \ref{thm:2.1}).

Another algorithm based on the Horner's method can be established using a modified argument in the computation of the classical Stirling numbers of the second kind shown in \cite{HS10}.  More precisely, we have the following algorithm. 

\begin{algorithm}\label{alg:1.2}
First, we may write the generalized Stirling numbers $S(n,k)=S(n,k,\alpha, \beta, r)$ defined by (\ref{eq:1.1**}) (see Definition \ref{def:1.0}) as 

\bn\label{eq:1.1***}
&&\la z\ra_{n,-\alpha}=\sum^n_{k=0} S(n,k)\la z-r\ra_{k,-\beta}\nonumber\\
&=&S(n,0)+(z-r)(S(n,1)+(z-r-\beta)(S(n,2)+(z-r-2\beta) (S(n,3)+\cdots \nonumber\\
&&(z-r-(n-1)\beta) S(n,n)))).
\en
Secondly, Use synthetic division to obtain $\la z\ra_{n,-\alpha}/(z-r)$, a polynomial of degree $\leq n-1$, with the remainder $S(n,0)$.
Then, evaluate $(\la z\ra_{n,-\alpha}/(z-r)-S(n,0) )/(z-r-\beta)$ to find the quotient polynomial of degree $\leq n-2$ 
as well as the remainder $S(n,1)$. Continue this process until a polynomial of degree $\leq 1$ 
left, which is $ S(n,n-1)+ (z-r-(n-1)\beta) S(n,n)$. A equivalent description of the above process can be presented as follows. 
Use Horner's method to find 

\[
f(r)\equiv \la z\ra_{n,-\alpha}=S(n,0)+(z-r) f_1(z),\quad deg \,\,f_1(z)\leq n-1,
\]
where the remainder is $S(n,0))$. Then, use Horner's method again to evaluate 

\[
f_1(z) =S(n,1)+(z-r-\beta) f_2(z), \quad deq\,\, f_2(z)\leq d-2,
\]
which generates the remainder $S(n,1)$. Continue the process and finally obtain 

\[
f_{n-1} =
S(n,n-1)+ (z-r-(n-1)\beta) S(n,n). 
\]

In short, we obtain $S(n,0)=\left. \la z\ra_{n,-\alpha}\right|_{z=r}$, $S(n,1)=\left. (\la z\ra_{n,-\alpha}-S(n,0))/(z-r)\right|_{z=r+\beta}$, etc.
\end{algorithm}
Algorithm \ref{alg:1.2} can be demonstrated by the following examples. 

\medbreak
\noindent{\bf Example 2.4} For the classical Stirling numbers of the second kind in the case of $n=5$ and $(\alpha, \beta, r)=(0,1,0)$, from expansion (\ref{eq:1.1***}) we have  
 
 \[
 z^5= S(5,0)+z(S(5,1) + (z-1) (S(5,2)+ (z-2) (S(5,3)+ (z-3) (S(5,4)+(z-4)S(5,5))))),
 \]
which implies $S(5,0)=0$ and 

\[
z^4= S(5,1) + (z-1) (S(4,2)+ (z-2) (S(4,3)+ (z-3) (S(4,4)+(z-4)S(5,5)))).
 \]
Thus, we may use the following division to evaluate $S(5,k)$ ($k=1,2,3,4,5$). 
 
 \begin{center}
\begin{tabular}{cccccc}
$1 \hspace{.23in}\vline$ & $1$ & $0$& $0$ & $0$ & $0$\\
 \rule[-3mm]{0mm}{8mm}
\hspace{.21 in} \vline&  & $1$& $1$& $1$ & $ 1$\\ 
 \hline \rule[-3mm]{0mm}{8mm}
$2\hspace{.19 in} \vline$ & $1$ & $1$ &$1$ & $1$ & {\bf 1}\\
  \rule[-3mm]{0mm}{8mm}
\hspace{.21 in} \vline&  & $2$& $6$& $14$ &\\ 
 \hline \rule[-3mm]{0mm}{8mm}
 $3 \hspace{.19in}\vline$& $1$ & $3$ &$7$ & {\bf 15} &\\
 \rule[-3mm]{0mm}{8mm}
\hspace{.21 in} \vline&  & $3$& $18$& & \\ 
 \hline \rule[-3mm]{0mm}{8mm}
 $4 \hspace{.19in}\vline$ & $1$ &$6$ & {\bf 25}& & \\
 \rule[-3mm]{0mm}{8mm}
 \hspace{.25in}\vline & & $4$ & & &\\
 \hline\rule[-3mm]{0mm}{8mm}
 &{\bf 1} &{\bf 10} &  & &\\
 \end{tabular}
 \end{center}
 \vspace{.15 in}
Hence, $S(5,1)=1$, $S(5,2)=15$, $S(5,3)=25$, $S(5,4)=10$, and $S(5,5)=1$. 

From (\ref{eq:1.1***}) we also immediately know that $S(n,n)=1$ because it is the coefficient of $z^n$ on the right-hand side while the coefficient 
on the left-hand side is $1$.

\medbreak

Let $\{ t_j\}^n_{j=1}$ be a strictly increasing $n$-sequence, and let $\sigma=\{ \sigma(j)\}^k_{j=1}$ be any strictly increasing integer sequence in $[1,n]$. There holds the following well-known refinement formula of divided difference (see, for example, \cite{deB}) 

\[
f[t, t-t_{\sigma(1)}, \ldots, t-t_{\sigma(k)}]=\sum^{\sigma(k)-k)}_{j=\sigma(1)-1}c(j)f[t,t_{j+1},\ldots, t-t_{j+k}],
\]
where $c(j)=c_{t,\sigma}>0$. Using this refinement formula one may obtain the refinement formula of the generalized Stirling numbers defined by (\ref{eq:1.1**}).

\begin{proposition}\label{pro:2.3}
Let $n,k\in{\bN}_0$ and the parameter triple $(\alpha, \beta, r)\not= (0,0,0)$ is in ${\bR}^3$ or ${\bC}^3$. Then there holds 
refinement formula, 

\[
\underline\triangle (\beta_{\sigma(1:k)}) \left. \la z\ra_{n,-\alpha}\right|_{z=r}=\sum^{\sigma(k)-k)}_{j=\sigma(1)-1}c(j) \underline\triangle (\beta_{j+1:j+k}) \left. \la z\ra_{n,-\alpha}\right|_{z=r},
\]
where 

\[
\underline\triangle (\beta_{\ell:j}) f:=f[t, t+\ell \beta, t+(\ell +1)\beta, \ldots, t+j\beta]
\]
\end{proposition}
\section{Generalized Stirling functions}
\setcounter{equation}{0}

We now extend the 
Stirling numbers $S(n,k))$ expressed by (\ref{eq:1.15-4}) to a more wider generation form using the idea of \cite{BKT03}. First, in order to cover as large a function class as possible, we recall that the {\it generalized fractional difference operator} $\Delta_\beta^{\eta, \ep}$ with an exponential factor, which is introduced in \cite{BKT03}.  More precisely, for $\eta\in {\bC}$, $\beta\in{\bR}_+$, $\epsilon \geq 0$, the generalized fractional difference operator $\Delta_\beta^{\eta, \ep}$ is defined for ``sufficient good'' functions $f$ by 

\be\label{eq:1.18}
\Delta^{\eta, \ep}_\beta f(z):=\sum_{j\geq 0} (-1)^j {\eta \choose j} e^{(\eta-j)\ep} f(z+(\eta -j)\beta)\quad (z\in {\bC}),
\ee
where ${\eta\choose j}$ are the general binomial coefficients given by 

\be\label{eq:1.19}
{\eta \choose j}=\frac{[\eta]_j}{j!}:=\frac{\eta (\eta-1)\cdots (\eta -j+1)}{j!}\quad (j\in {\bN}),
\ee
with $[\beta]_0=1$. Noting the generalized Stirling numbers $S(n,k)$ can be represented by (\ref{eq:1.15}), 
or equivalently, 
\[
S(n,k)=\frac{1}{\beta^kk!} \lim_{z\to r} \Delta^k_\beta \la z\ra_{n,-\alpha},
\]
which has an extension shown in (\ref{eq:1.16**}). We now extend (\ref{eq:1.16**}) to a more generalized form shown in the following definition.

\begin{definition}\label{def:2.2}
The generalized 
Stirling functions, $S(\ga, \eta, $ $\alpha,\beta,r; \ep)$ for any complex numbers $\ga$ and $\eta$ are given by 

\be\label{eq:1.20}
S(\ga,\eta;\ep)\equiv S(\ga,\eta,\alpha,\beta,r; \epsilon):=\frac{1}{\beta^\eta \Ga(\eta +1)}\lim_{z\to r} \Delta^{\eta, \ep}_\beta (\langle z\rangle_{\ga,-\alpha})\quad (\ep \geq 0),
\ee
provided the limit exists; or equivalently, by 

\be\label{eq:1.21}
S(\ga,\eta;\ep )\equiv S(\ga,\eta,\alpha,\beta,r;\ep )=\frac{1}{\beta^\eta \Ga(\eta+1)}\sum_{j\geq 0} (-1)^j {\eta\choose j}e^{(\eta-j)\ep}\langle r+(\eta-j)\beta \rangle_{\ga,-\alpha}\quad (\ga\not= 0),
\ee
provided the series converges absolutely. 
and 

\be\label{eq:1.21*}
S(0,\eta)=\frac{(e^\epsilon -1)^\eta}{\beta^\eta \Gamma (\eta +1)}.
\ee
\end{definition}
From (\ref{eq:1.21}), we immediately have 

\be\label{eq:1.21**}
S(\ga,0;\epsilon)=\la r\ra_{\ga,-\alpha} \quad (\ga\not= 0).
\ee

Now, an explicit expression of $S(\ga,\eta;\ep)$ can be given by the following result.

\begin{theorem}\label{thm:4.8}
If $\ga \in{\bC}$ and either of the conditions $\eta\in{\bC}$ ($\eta\notin {\bZ}$), $\ep>0$, or $\eta\in{\bC}$ ($\eta\notin {\bZ},$ $Re (\eta)>Re (\ga)$), $\ep=0$ hold, then the generalized Stirling functions $S(\ga, \eta;\ep)$ can be represented in the form (\ref{eq:1.21})
and $S(0,\eta;\ep)=\delta_{\eta,0}$. In particular, if $\ga=n\in {\bN}_0$, $\eta =k\in {\bN}$, and $\ep \geq 0$, then the corresponding generalized Stirling functions $S(n,k;\ep)$ has the representation (\ref{eq:1.21}). 
\end{theorem}

We now present the recurrence relation of the generalized Stirling functions defined by (\ref{eq:1.21}) by using the recurrence relations of  the generalized raising and falling factorial functions shown in Proposition \ref{pro:2.2}. 

\begin{theorem}\label{thm:4.9}
There hold the following three results.

(a) For $\ga \in{\bC}$, $\eta\in{\bC}$ ($\eta\notin {\bZ}$), and $\ep>0$, the generalized 
Stirling functions $S(\ga, \eta;\ep)$ defined by (\ref{eq:1.21}) satisfy 

\be\label{eq:1.23}
S(\ga, \eta;\ep)=(r+\eta \beta -(\ga -1)\alpha)  S(\ga-1,\eta;\ep)+S(\ga-1,\eta -1;\ep).
\ee

(b) Let $\ga\in{\bC}$, $\eta\in{\bC}$ ($\eta\notin {\bZ}),$ and $Re (\eta)>Re (\ga)$). The generalized Stirling functions $S(\ga,\eta)\equiv S(\ga,\eta;0)$ satisfy  

\be\label{eq:1.24}
S(\ga,\eta)=(r+\eta \beta-(\ga-1)\alpha)S(\ga-1,\eta)+S(\ga-1,\eta-1).
\ee

(c) For $\ga \in{\bC}$, $k\in{\bN}$, and $\ep\geq 0$, the generalized Stirling functions $S(\ga, k;\ep;h)$ defined by (\ref{eq:1.21}) satisfy 

\be\label{eq:1.25}
S(\ga, k;\ep)=(r+k \beta -(\ga -1)\alpha)  S(\ga-1,k;\ep)+S(\ga-1,k -1;\ep).
\ee
In particular, 

\[
S(\ga, k)=(r+k \beta -(\ga -1)\alpha)  S(\ga-1,k)+S(\ga-1,k -1).
\]
\end{theorem}

Clearly, Theorem 6 in \cite{BKT03} is a special case of Theorem \ref{thm:4.9} for $\alpha, \beta=0$. And Theorem 3, Corollaries 3.1 and 3.2 in \cite{BKT03} are special cases of Theorem \ref{thm:4.9} for $\alpha, \beta=0$ and $\ga =n\in{\bN}$. 

Now we construct the exponential generating function for the generalized Stirling functions $S(n, \eta;\ep)$.

\begin{theorem}\label{thm:4.10}
Let $z \in{\bC}$, $\eta\in{\bC}$, and $\ep\geq 0$. The generating function for the generalized Stirling functions $S(\ga, \eta;\ep)$ defined by (\ref{eq:1.21}) with $\ga=n$ and $\alpha \beta \not= 0$ is 

\be\label{eq:3.1}
\frac{1}{\Gamma (\eta +1)}(1+\alpha z)^{r/\alpha} \left( \frac{e^\epsilon(1+\alpha z)^{\beta/\alpha} -1}{\beta}\right)^\eta = \sum_{n\geq 0} S(n,\eta; \epsilon) \frac{z^n}{n!}
\ee
for $\eta \not\in {\bZ}$ and $\epsilon >0$, and 

\be\label{eq:3.2}
\frac{1}{k!}(1+\alpha z)^{r/\alpha} \left( \frac{e^\epsilon(1+\alpha z)^{\beta/\alpha} -1}{\beta}\right)^k = \sum_{n\geq 0} S(n,k; \epsilon) \frac{z^n}{n!}
\ee
for $\eta =k\in {\bN}_0$ and $\epsilon \geq 0$.
\end{theorem}

\noindent{\bf Remark 3.1} The condition $\alpha \beta\not= 0$ is not necessary for the left-hand side of (\ref{eq:3.1}). In fact, taking $r=0$, $\beta =1$, and letting $\alpha \to 0^+$, we see that (\ref{eq:3.1}) yields the generating function for the generalized Stirling functions of the second kind: 

\[
\frac{1}{\Ga(\eta+1)} (e^{z+s}-1)^\eta =\sum _{n\geq 0} S(n,\eta, 0, 1,0;\ep) \frac{z^n}{n!},
\]
which was studied in Theorem $4$ 
of \cite{BKT03}, and it can be considered as a particular case of our Theorem \ref{thm:4.10}. 

Similarly, taking $\ep, r=0$, $\alpha =1$ and letting $\beta\to 0^+$ yields the generating function of the generalized  Stirling functions of the first kind:

\[
\frac{1}{\Ga (\eta +1)} (\ln (1+z))^\eta =\sum_{n\geq 0} S(n,\eta, 1,0,0)\frac{z^n}{n!}.
\]

\section{More properties of the generalized Stirling functions and numbers}
\setcounter{equation}{0}

let us consider the set of formal power series (f.p.s.) $\F = {\bR}[[t; \{c_k\}]]$ or ${\bC}[[t;\{ c\}]]$ (where $c=(c_0, c_1, c_2, \ldots)$ satisfies $c_0=1$, $c_k>0$ for all $k=1,2,\ldots$); the \emph{order} of $f(t) \in \F$, $f(t)=\sum_{k=0}^\infty f_kt^k/ c_k$, is the minimal number $r\in {\bN}$ such that $f_r \neq 0$; $\F_r$ is the set of formal power series of order $r$. It is known that $\F_0$ is the set of \emph{invertible} f.p.s. and $\F_1$ is the set of \emph{compositionally invertible} f.p.s., that is, the f.p.s.'s $f(t)$ for which the compositional inverse $\overline{f}(t)$ exists such that $f(\overline{f}(t)) = \overline{f}(f(t)) = t$. We call the element $g\in \F$ with the form $g(x)=\sum_{k\geq 0} \frac{x^k}{c_k}$ a generalized power series (GPS) associated with $\{ c_n\}$ or, simply, a (c)-GPS, and $\F$ the GPS set associated with $\{ c_n\}$. In particular, when $c=(1,1, \ldots)$, the corresponding $\F$ and $\F_r$ denote the classical formal power series and the classical formal power series of order $r$, respectively.

We now develop a kind of asymptotic expansions for the generalized Stirling functions $S(n,\mu, r;\ep)\equiv S(n, \mu, \alpha, \beta, r; \ep)$ and $S(n,\mu, \mu r;\ep)\equiv S(n, \mu, \alpha, \beta,\mu r; \ep)$ and generalized Stirling numbers $S(n+\mu, \mu, r)\equiv S(n+\mu, \mu, \alpha,\beta, r)$ and $S(n+\mu, \mu, \mu r)\equiv S(n+\mu, \mu, \alpha,\beta, \mu r)$ for large $\mu$ and $n$ with the condition $n=0(\mu^{1/2})$ ($\mu\to \infty$). The  asymptotic expansions of Hsu and Shiue Stirling numbers in \cite{HS98} and Tsylova Stirling numbers in \cite{Tsy}, involving a generalization of Moser and Wyman's result \cite{MW}, are included as particular cases. 

The major tool of construction of the asymptotic expansion is the known result about the asymptotic formula for the coefficients of power-type generating functions involving large parameters shown in \cite{Hsu90}. Let $\sigma(n)$ be the set of partition of $n$ ($n\in{\bN})$, which can be represented by $1^{k_1}2^{k_2}\cdots n^{k_n}$ with $1k_1+2k_2+\cdots nk_n=n$, $k_j\geq 0$ ($j=1,2,\ldots, n$), and with $k=k_1+k_2+\cdots +k_n$ expressing the number of the parts of the partition. For given $k$ ($1\leq k\leq n$), we denote by $\sigma(n,k)$ the subset of $\sigma(n)$ consisting of partitions of $n$ having $k$ parts. 

Let $\phi(z)=\sum_{n\geq 0} a_n z^n$ be a formal power series over the complex field ${\bC}$ in $\F_0$, with $a_0=g(0)=1$. For every $j$ ($0\leq j< n$) define 

\be\label{eq:5.1}
W(n,j)=\sum_{\sigma (n,n-j)} \frac{a_1^{k_1} a_2^{k_2}\cdots a_n^{k_n}}{k_1! k_2!\cdots k_n!},
\ee
where the summation is taken over all such partition $1^{k_1} 2^{k_2}\cdots n^{k_n}$ of $n$ that have $n-j$ parts. 
We have the following known result (see for instance \cite{Hor85}): 

For a fixed $m\in {\bN}$ and for large $\mu$ and $n$ such that $n=o(\mu^{1/2})$ ($\mu \to \infty$), we have the asymptotic expansion 

\be\label{eq:5.2}
\frac{1}{[\mu]_n} [z^n] (\phi(z))^\mu =\sum^m_{j=0} \frac{ W(n,j)}{[\mu-n+j]_j}+o\left( \frac{ W(n,m)}{[\mu-n+m]_m}\right),
\ee
where $W(n,j)$ are given by (\ref{eq:5.1}). (\ref{eq:5.2}) is used to derive the Hsu-Shiue Stirling numbers in \cite{HS98}. We now generalize (\ref{eq:5.2}) and the corresponding argument to give asymptotic expansion formulas of generalized Stirling functions $S(n,\mu, r;\ep)\equiv S(n, \mu, \alpha, \beta,$ $r; \ep)$, $S(n,\mu, \mu r;\ep)\equiv S(n, \mu, \alpha, \beta,\mu r; \ep)$, $S(n+\mu, \mu, r)\equiv S(n+\mu, \mu, \alpha,\beta, r)$ and $S(n+\mu, \mu, \mu r)\equiv S(n+\mu, \mu, \alpha,\beta, \mu r)$ for large $\mu$ and $n$ with the condition $n=0(\mu^{1/2})$ as $\mu\to \infty$. 

Let $g(z)=\sum_{n\geq 0} a_n z^n$ be a formal power series over the complex field ${\bC}$ in $\F_0$, with $a_0=g(0)\not= 0$. We may write 

\[
g(z)=a_0 \sum_{n\geq 0} \frac{a_n}{a_0} z^n.
\]

For a fixed $m\in {\bN}$ and for large $\mu$ and $n$ such that $n=o(\mu^{1/2})$ ($\mu \to \infty$), From  formulas (\ref{eq:5.1}) and (\ref{eq:5.2}) we have the asymptotic expansion 

\be\label{eq:5.2*}
\frac{1}{[\mu]_n} [z^n] (g(z))^\mu =\sum^m_{j=0} \frac{W(n,j)}{a_0^{n-\mu-j}[\mu-n+j]_j}+o\left( \frac{W(n,m)}{a_0^{m-\mu}[\mu-n+m]_m}\right),
\ee
where $W(n,j)$ are given by (\ref{eq:5.1}). 
In particular, when $n$ is fixed, the remainder estimate becomes $O(\mu^{-m-1})$. 

To apply (\ref{eq:5.2}) to the generalized Stirling numbers $S(\ga, \eta;\ep)$ defined by (\ref{eq:1.21}) with $\ga=n$, $\eta=\mu$ and $\alpha \beta \not= 0$, let us use (\ref{eq:3.2}) to take 

\be\label{eq:5.3}
g(z)= (1+\alpha z)^{r/\alpha}\frac{e^\ep (1+\alpha z)^{\beta/\alpha} -1}{\beta } =\sum_{n\geq 0} \frac{S(n, 1;\ep)}{n!} z^n
\ee
when $\ep\not= 0$, and 

\be\label{eq:5.3*}
\bar g(z)= (1+\alpha z)^{r/\alpha}\frac{(1+\alpha z)^{\beta/\alpha} -1}{\beta z} =\sum_{n\geq 0} \frac{S(n+1, 1)}{(n+1)!} z^n
\ee
when $\ep=0$, so that $g(0)=(e^\ep-1)/\beta$ ($\ep\not= 0$) and $\bar g(0)=1$ ($\ep =0$) not being zero in both cases, where $S(n,1;\ep)\equiv S(n,1,\alpha, \beta, r;\ep)$ and $S(n+1,1)\equiv S(n+1,1,\alpha,\beta,r)$, $g(0)=(e^\ep-1)/\beta$. Consequently, from (\ref{eq:3.2}) we have 

\bn\label{eq:5.4}
(g(z))^\mu &=&(1+\alpha z)^{\mu r/\alpha} \left( \frac{e^\ep (1+\alpha z)^{\beta/\alpha}-1}{\beta }\right)^\mu\nonumber \\&=& 
\mu! \sum_{n\geq 0} \frac{S(n,\mu,\alpha,\beta,\mu r;\ep)}{n!} z^n
\en
for $\ep\not= 0$, and 

\bn\label{eq:5.4*}
 (\bar g(z))^\mu &=&(1+\alpha z)^{\mu r/\alpha} \left( \frac{(1+\alpha z)^{\beta/\alpha}-1}{\beta z}\right)^\mu\nonumber \\
&=& \mu! \sum_{n\geq 0} \frac{S(n+\mu,\mu,\alpha,\beta,\mu r)}{(n+\mu)!} z^n
\en
for $\ep=0$. Therefore, making use of (\ref{eq:5.2*}) yields

\bn\label{eq:5.5}
&&\frac{S(n,\mu,\alpha,\beta,\mu r;\ep)}{[\mu]_n[n]_\mu}\nonumber\\
&=&\left( \frac{\beta}{e^\ep -1}\right)^{n-\mu}\sum^m_{j=0} \left( \frac{e^\ep -1}{\beta}\right)^j\frac{W(n,j)}{[\mu-n+j]_j}+o\left( \left(\frac{\beta}{e^\ep -1}\right)^{n-\mu}\frac{W(n,m)}{[\mu -n+m]_m}\right)
\nonumber\\
\en
for $\ep\not= 0$, and 

\be\label{eq:5.5*}
\frac{S(n+\mu,\mu,\alpha,\beta,\mu r)}{[\mu]_n[n+\mu]_\mu}= \sum^m_{j=0} \frac{W(n,j)}{[\mu-n+j]_j}+o\left( \frac{W(n,m)}{[\mu -n+m]_m}\right)
\ee
for $\ep=0$, where $n=o(\mu^{1/2})$ as $\mu\to\infty$ and $W(n,j)$ ($j=0,1,2,\ldots$) are given by (\ref{eq:5.1}) with $a_j$ being determined by (\ref{eq:5.3}); namely, for $\ep\not= 0$, $a_0=(e^\ep-1)/\beta$ and 

\be\label{eq:5.6}
a_j=[z^j]g(z)=\frac{S(j,1;\ep)}{j!},
\ee
while for $\ep=0$, $a_0=1$ and 

\be\label{eq:5.6*}
a_j=[z^j] \bar g(z)=\frac{S(j+1,1)}{(j+1)!}.
\ee
The coefficients defined by (\ref{eq:5.6}) and (\ref{eq:5.6*}) can be evaluated by using the Vandermonde-Chu formula as follows. From (\ref{eq:5.3}), for $j=1,2,\ldots,$ we have 

\bns
&&[z^j]g(z)=[z^j](1+\alpha z)^{r/\alpha} \left[ \frac{e^\ep -1}{\beta }+\frac{e^\ep}{\beta }\sum_{k\geq 1} {\beta/\alpha \choose k} (\alpha z)^k \right]\\
&=&[z^j] \left[ \frac{e^\ep -1}{\beta} \sum_{\ell \geq 0} {r/\alpha \choose \ell} (\alpha z)^\ell +
\frac{e^\ep}{\beta} \sum_{\ell \geq 0} \sum_{k \geq 1} {r/\alpha \choose \ell}{\beta/\alpha \choose k}(\alpha z)^{\ell+k}\right]\\
&=&\frac{e^\ep -1}{\beta} \alpha^j {r/\alpha \choose j}+\frac{e^\ep}{\beta}\alpha^j \sum^j_{k=1}{r/\alpha \choose {j-k}}{\beta/\alpha \choose k}\\
&=&\frac{e^\ep -1}{ j!\beta }\la r\ra_{j, -\alpha}+\frac{e^\ep}{\beta}\alpha^j \left[ {r/\alpha +\beta/\alpha \choose j}-{r/\alpha \choose j}\right]\\
&=&\frac{1}{j!\beta} \left[\la r+\beta\ra_{j,-\alpha}+(e^\ep -2)\la r\ra_{j,-\alpha}\right].
\ens
Here, the classical Vandermonde-Chu convolution formula we used above, regarded as ``perhaps the most widely used combinatorial identity'' (see P. 8 in \cite{Rio79} by Riordan and PP. 51, 61, 64, and 227 in \cite{And} by Andrews), which can be written as 

\[
\sum^n_{k=0} {x\choose k}{y\choose {n-k}}={x+y\choose n}\quad (x,y\in{\bR}, n\in {\bN}_0).
\]
Similarly, we obtain

\[
[z^j] \bar g(z)= \frac{1}{(j+1)!\beta} \left[ \la r+\beta\ra_{j+1, -\alpha}-\la r\ra_{j+1, -\alpha}\right]
\]
for $j=0,1,2,\ldots.$ Hence, we may survey the above into the following theorem.

\begin{theorem}\label{thm:5.1}
For $\ep\not= 0$, there holds the asymptotic expansion (\ref{eq:5.5}) of $S(n,\mu, \mu r;\ep)\equiv S(n, \mu, \alpha, \beta,\mu r; \ep)$ for $n$ with $n=o(\mu^{1/2})$ ($\mu\to \infty$), where $W(n,j)$ is defined by (\ref{eq:5.1}) with $a_0=(e^\ep -1)/\beta$ and 

\[
a_j=\frac{1}{j!\beta} \left[\la r+\beta\ra_{j,-\alpha}+(e^\ep -2)\la r\ra_{j,-\alpha}\right]\quad (j=1,2,\ldots).
\]
For $\ep=0$, there holds the asymptotic expansion (\ref{eq:5.5*}) of $S(n+\mu,\mu, \mu r)\equiv S(n+\mu, \mu, \alpha, \beta,\mu r)$ for $n$ with $n=o(\mu^{1/2})$ ($\mu\to \infty$), where $W(n,j)$ is defined by (\ref{eq:5.1}) with

\[
a_j=\frac{1}{(j+1)!\beta} \left[ \la r+\beta\ra_{j+1, -\alpha}-\la r\ra_{j+1, -\alpha}\right]\quad j=0,1,\ldots.
\]
\end{theorem}

Since the formulas (\ref{eq:5.5}) and (\ref{eq:5.5*}) with $W(n,j)$ and $a_j$ presented in (\ref{eq:5.1}) and Theorem \ref{thm:5.1}, respectively, are algebraic analytic identities, we may replace $r$ by $r/\mu$ in the formulas and obtain the following corollary.

\begin{corollary}
For $\ep\not= 0$, by replacing the quantity $r$ by $r/\mu$, the asymptotic expansion (\ref{eq:5.5}) is also applicable to $S(n,\mu, r;\ep)\equiv S(n, \mu, \alpha, \beta,r; \ep)$ for $n$ with $n=o(\mu^{1/2})$ ($\mu\to \infty$), where $W(n,j)$ is defined by (\ref{eq:5.1}) with $a_0=(e^\ep -1)/\beta$ and 

\[
a_j=\frac{1}{j!\beta} \left[\left \la \frac{r}{\mu}+\beta\right \ra_{j,-\alpha}+(e^\ep -2)\left \la \frac{r}{\mu}\right\ra_{j,-\alpha}\right]\quad (j=1,2,\ldots).
\]
For $\ep=0$, by replacing the quantity $r$ by $r/\mu$, the asymptotic expansion (\ref{eq:5.5*}) is also applicable to   $S(n+\mu,\mu,  r)\equiv S(n+\mu, \mu, \alpha, \beta,r)$ for $n$ with $n=o(\mu^{1/2})$ ($\mu\to \infty$), where $W(n,j)$ is defined by (\ref{eq:5.1}) with

\[
a_j=\frac{1}{(j+1)!\beta} \left[ \left \la \frac{r}{\mu}+\beta\right \ra_{j+1, -\alpha}-\left \la \frac{r}{\mu}\right\ra_{j+1, -\alpha}\right]\quad j=0,1,\ldots.
\]
\end{corollary}

In the recent literature, special emphasis has been given to the concept of {\em Riordan arrays},
which are a generalization of the well-known Pascal triangle. Riordan arrays are infinite, lower
triangular matrices defined by the generating function of their columns. They form a group, called
{\em the Riordan group} (see Shapiro et al. \cite{SGWW}). Some of the main results on the Riordan
group and its application to combinatorial sums and identities can be found in Sprugnoli
\cite{Spr1, Spr2}, on subgroups of the Riordan group in Peart and Woan \cite{PW} and Shapiro
\cite{Sha0}, on some characterizations of Riordan matrices in Rogers \cite{Rog}, Merlini et al.
\cite{MRSV} and He et al. \cite{HS09}, and on many interesting related results in Cheon et al. \cite{CKS, CKS2}, He et al. \cite{HHS}, Nkwanta \cite{Nkw}, Shapiro \cite{Sha1, Sha2}, and so forth. We now generalize the Riordan arrays associated with classical power series to those associated with (c)-GPS, where $c=\{ c_k=k!\}_{k\geq 0}$. The Riordan arrays associated with other (c)-GPS can be found in author's later paper. More precisely, let $c=\{ c_k=k!\}_{k\geq 0}$. The (c)-Riordan array generated by $d(t)\in \F_0$ and $h(t)\in \F_1$ with respect to $\{ c_k\}_{k\geq 0}$ is an infinite complex matrix $[d_{n,k}]_{0\leq k\leq n}$, whose bivariate generating function has the form 

\begin{equation}\label{eq:def3}
F(t,x)=\sum_{n,k} d_{n,k} \frac{t^n}{n!} x^k =d(t) e^{xh(t)},
\end{equation}
which is called a {\it Sheffer type Riordan array}.

Thus, the $(n,k)$ entry of (c)-Riordan array $[d_{n,k}]$ is 

\begin{equation}\label{eq:dnk}
d_{n,k}=\left[ \frac{t^n}{n!}\right] d(t) \frac{(h(t))^k}{k!}=[t^n]\frac{n!}{k!} d(t) (h(t))^k
\end{equation}
for all $0\leq k\leq n$ and $d_{n,k}=0$ otherwise. It is easy to see that a lower triangular array $[d_{n,k}]$ is a (c)-Riordan array if and only if the array $(k! d_{n,k}/n!)$ is
a (1)-Riordan array, i.e., a classical Riordan array. Evidently, $[d_{n,k}]=(d(t), h(t))$ can be written as 

\begin{equation}\label{eq:RM}
[d_{n,k}]=D [ [t^n] d(t) (h(t))^k]_{n\geq k\geq 0}D^{-1},
\end{equation}
where $D=diag (1, 1, 2!,\ldots)$.

Rogers \cite{Rog} introduced the concept of the $A$-{\em sequence} for
the classical Riordan arrays; Merlini et al. \cite{MRSV} introduced the related concept of the $Z$-{\em sequence}
and showed that these two concepts, together with the element $d_{0,0}$, completely characterize a
proper classical Riordan array. In \cite{HS09}, Sprugnoli and the author consider the characterization of Riordan arrays, their multiplications, and their inverses by means of the $A$- and $Z$-sequences.

In \cite{Rog}, Rogers states that for every proper Riordan array $D = (d(t),\,h(t))$
there exists a sequence $A=(a_k)_{k\in{\bN}}$ such that for every $n,k\in{\bN}$ we have:
  \bn\label{AseqRA}
&&[t^{n+1}]d(t) (h(t))^{k+1} \nonumber\\
&=& a_0[t^n]d(t) (h(t))^{k} + a_1[t^n]d(t) (h(t))^{k+1}+ a_2[t^n]d(t) (h(t))^{k+2} + \cdots \nonumber\\
&=& \sum_{j=0}^\infty
  a_j[t^n]d(t) (h(t))^{k+j}
\en
where the sum is actually finite since $d_{n,k}=0,\ \forall k>n$. We can reformulate it to the generalized (c)-Riordan array as follows.

\begin{theorem}\label{thm:Aseq}
An infinite lower triangular array $D=\matrice{d}=(d(t), h(t))$ is a (c)-Riordan array if and only if a sequence $A
= (a_0\not= 0, a_1, a_2, \ldots)$ exists such that for every $n,k\in{\bN}$ relation 
\begin{equation}\label{AseqGRA}
\frac{c_{k+1}}{c_{n+1}}d_{n+1,k+1}=\frac{c_0}{c_n}a_0d_{n,k} + \frac{c_1}{c_n}a_1d_{n,k+1} + \frac{c_2}{c_n}a_2d_{n,k+2} + \cdots = \sum_{j=0}^\infty \frac{c_{k+j}}{c_n}  a_jd_{n,k+j}
\end{equation}
holds. In addition, the generating function $A(t)$ of $A-$ sequence is uniquely determined by $t A(h(t))=h(t)$.
\end{theorem}

We now use Theorem \ref{thm:Aseq} to establish a new recursive relationship of generalized Stirling numbers. From expression (\ref{eq:3.2}) in Theorem \ref{thm:4.10} with $\ep=0$ and $\alpha\beta\not= 0$, we have the generating function of the generalized Stirling numbers shown below:

\be\label{eq:5.7}
\frac{1}{k!}(1+\alpha z)^{r/\alpha} \left( \frac{(1+\alpha z)^{\beta/\alpha} -1}{\beta}\right)^k = \sum_{n\geq 0} S(n,k) \frac{z^n}{n!}. 
\ee

\begin{theorem}\label{SRA}
Let $\alpha\beta\not= 0$. The $A-$ sequence $(a_n)_{n\in{\bN}_0}$ of the Riordan array of the generalized Stirling number array $[d_{n,k}=k! S(n,k)/n!]_{0\leq k\leq n}$ satisfies 

\be\label{eq:5.7*}
a_0 =1,\quad a_n=-\frac{1}{\alpha} \sum^n_{k=1} a_{n-k} \frac{\la \alpha \ra_{k+1, -\beta}}{(k+1)!}
\ee
for all $n\geq 1$.
\end{theorem}

To find the fist column of the array $[d_{n,k}]_{0\leq k\leq n}$, we consider (\ref{eq:5.7}) for $k=0$ and have 

\[
(1+\alpha z)^{r/\alpha} =\sum_{n\geq 0} \frac{S(n,0)}{n!} z^n.
\]
On the other hand, 

\[
(1+\alpha z)^{r/\alpha}=\sum_{n\geq 0} {r/\alpha \choose n} (\alpha z)^n.
\]
Comparing the right-hand sides of the last two equations, we obtain 

\be\label{eq:5.8*}
S(n,0)\equiv S(n,0,\alpha, \beta, r)=n! {r/\alpha \choose n} \alpha^a =\la r\ra_{n,-\alpha}.
\ee
Formula (\ref{eq:5.8*}) was given in (\ref{eq:1.5}) and also in (\ref{eq:1.16**}), which are derived by different approaches.

From (\ref{eq:5.7}) we have 

\be\label{eq:5.8**}
[d_{n,k}]_{0\leq k\leq n}=\left[ \frac{k!}{n!} S(n,k)\right]_{0\leq k\leq n},
\ee
where $S(n,k)\equiv S(n,k,\alpha, \beta, r)$ ($\alpha\beta\not= 0$). Therefore, surveying the above process, we obtain an algorithm to evaluate generalized Stirling numbers $S(n,k)\equiv S(n,k,\alpha, \beta, r)$ with $\alpha\beta\not= 0$. 

\begin{algorithm}\label{alg:4.1}
Denote $d(t)= (1+\alpha z)^{r/\alpha}$ and $h(z)=((1+\alpha z)^{\beta/\alpha}-1)/\beta$ ($\alpha\beta\not= 0$). Let $n,k\in{\bN}_0$ and $\alpha\beta\not= 0$. Then we may find $A$-sequence $(a_n)_{n \in {\bN}_0}$ shown in (\ref{eq:5.7*}) and establish the array (\ref{eq:5.8**}) except its first column by using the recursive relation (\ref{AseqGRA}) shown in Theorem \ref{thm:Aseq}, i.e., 

\be\label{eq:5.9}
\frac{k!}{n!} S(n,k)=\sum_{j\geq 0}a_j \frac{(k+j-1)!}{(n-1)!} S(n-1, k+j-1)
\ee
for all $1\leq k\leq n$. The first column of array (\ref{eq:5.8**}) can be constructed by using (\ref{eq:5.8*}). Thus, the $n$th entry of the first column is  

\be\label{eq:5.10}
\frac{1}{n!} S(n,0)=\frac{\la r\ra_{n,-\alpha}}{n!}.
\ee
Finally, all $S(n,k)\equiv S(n,k,\alpha,\beta,r)$ ($0\leq k\leq n$) can be read from a modification of array (\ref{eq:5.8**}); namely from 

\[
\left[ \frac{n!}{k!} d_{n,k}\right]_{0\leq k\leq n}=\left[ S(n,k)\right]_{0\leq k\leq n},
\]
where $S(n,k)=n\sum_{j\geq 0}a_j  [k+j-1]_{j-1} S(n-1, k+j-1)$ when $1\leq k\leq n$, 
and $S(n,0)$ can be obtained from (\ref{eq:5.10}) or (\ref{eq:5.8*}). 
\end{algorithm}
\medbreak

\noindent{\bf Remark 4.1} Similar to the argument in Remark 3.1, the condition $\alpha\beta \not= 0$ in Theorem \ref{SRA} and Algorithm \ref{alg:4.1} is not necessary. Algorithm \ref{alg:4.1} can be modified to adapt some of cases when $\alpha\beta =0$. We will show the application of Algorithm \ref{alg:4.1} to the calculations of the classical Stirling numbers of the second and the first kind, i.e., $S(n,k, \alpha, \beta, r))=S(n,k,0,1,0)$ and $S(n,k,\alpha, \beta, r)=S(n,k,1,0,0)$, in Examples 4.2 and 4.3, respectively. 

\medbreak

\noindent{\bf Example 4.1} 
For the Howard's weighted degenerated Stirling numbers $S(n,k)\equiv S(n,k, 1,1,-1)$. From Algorithm \ref{alg:4.1} or Theorem \ref{SRA}, we immediately have generating function of the corresponding $A$-sequence $A(z)=1$. 
Then, using (\ref{eq:5.9}) and (\ref{eq:5.10}) we obtain the Riordan array $[d_{n,k}]_{0\leq k\leq n}=\left[ \frac{k!}{n!} S(n,k)\right]_{0\leq k\leq n}$ as 

\[
\left[ \frac{k!}{n!} S(n,k)\right]_{0\leq k\leq n}=\left[ \begin{array}{rrrrr} 1 &  & & &  \\
-1 &1 & & & \\ 1 &-1& 1& &  \\ -1 & 1& -1 &1 & \\ 1 & -1& 1& -1& 1\end{array}\right].
\]
Therefore, 

\[
\left[ S(n,k)\right]_{0\leq k\leq n}=\left[\begin{array}{rrrrr} 1 &  & & &  \\
-1 &1 & & & \\ 2 &-2& 1& &  \\ -6 & 6& -3 &1 & \\ 24 & -24& 12& -4& 1\end{array}\right],
\]
which gives $S(0,0)=1$; $S(1,0)=-1$, $S(1,1)=1$; $S(2,0)=2$, $S(2,1)=-2$, $S(2,2)=1$; $S(3,0)=-6$, $ S(3,1)=6$, $S(3,2)=-3$, $S(3,3)=1$; and $S(4,0)=24$, $S(4,1)=-24$, $S(4,2)=12$, $S(4,3)=-4$, and $S(4,4)=1$ row by row. 
\medbreak

\noindent{\bf Example 4.2} As we have presented in Remarks 3.1 and 4.1, the condition $\alpha\beta \not= 0$ in Theorems \ref{thm:4.10} and \ref{SRA} and Algorithm \ref{alg:4.1} is not necessary. Here, we demonstrate how to modify Algorithm \ref{alg:4.1} for the case of $(\alpha, \beta ,r) =(0,1,0)$. The generating function of the corresponding classical Stirling numbers $\{ S(n,k)\equiv S(n,k,0,1,0)\}_{0\leq k\leq n}$ of the second kind is 

\[
\frac{1}{k!} (e^z-1)^k =\sum_{n\geq 0} S( n,k)\frac{z^n}{n!}.
\]
Thus the corresponding Riordan array has generating functions $d(z)=1$ and $h(z)=e^z-1$. Since the compositional inverse of $h(z)$ is $\bar h(z)=\ln (1+z)$, the $A$-sequence characterization of the Riordan array has generating function 

\[
A(z)=\frac{z}{\ln (1+z)}=\frac{z}{\sum_{k\geq 1} \frac{(-1)^{k-1}}{k} z^k}=\frac{1}{\sum_{k\geq 0} \frac{(-1)^k}{k+1} z^k},
\]
which coefficients $\{ a_n\}_{n\geq 0}$, i.e., the elements of $A$-sequence, can be solved from the above equation as 

\[
a_0=1, \quad a_n=-\sum^n_{k=1} a_{n-k} \frac{(-1)^k}{k+1}=\sum^{n+1}_{k=2} a_{n-k+1}\frac{(-1)^k}{k}\quad (n\geq 1).
\]
Thus, we obtain the first few $a_n$:

\[
a_0=1, \,\, a_1=\frac{1}{2}, \,\, a_2=-\frac{1}{12},\,\, a_3=\frac{1}{24},\,\, a_4=-\frac{19}{720}, \,\, etc.
\]
Similar to Algorithm \ref{alg:4.1}, we may find the Riordan array 

\[
[d_{n,k}]_{0\leq k\leq n}=\left[ \frac{k!}{n!} S(n,k)\right]_{0\leq k\leq n}=\left[ \frac{k!}{n!} S(n,k)\right]_{0\leq k\leq n}=\left[ \begin{array}{rrrrr} 1 &  & & &  \\
0 &1 & & & \\ 0 &\frac{1}{2}& 1& &  \\0 & \frac{1}{6}& 1 &1 & \\ 0 & \frac{1}{24}& \frac{7}{12}& \frac{3}{2}& 1\end{array}\right].
\]
The Riordan Stirling array of the Stirling numbers of the second kind is 

\[
\left[ S(n,k)\right]_{0\leq k \leq n}= \left[ \frac{k!}{n!} S(n,k)\right]_{0\leq k\leq n}=\left[ \begin{array}{rrrrr} 1 &  & & &  \\
0 &1 & & & \\ 0 &1& 1& &  \\ 0 & 1& 3 &1 & \\ 0 & 1& 7& 6& 1\end{array}\right],
\]
which gives all $S(n,k)=S(n,k,0,1,0)$ for $0\leq k \leq 4$. For instance, $S(4,0)=0$, $S(4,1)=1$, $S(4,2)=7$, $S(4,3)=6$, and $S(4,4)=1$.

{\bf Example 4.3} For $(\alpha, \beta, r)=(1,0,0)$, we can also applied a modification of Algorithm \ref{alg:4.1} to evaluate the classical Stirling numbers of the first kind $s(n,k)\equiv S(n,k,1,0,0)$ as follows.
In this case, we have the corresponding Riordan array $(d(z), h(z))=(1, \ln (1+z))$. Thus the compositional inverse of $\bar h(z)=e^z-1$. Thus the $A$-sequence $\{ a_n\}_{n\geq 0}$ has its generating function 

\[
A(z)=\frac{z}{\bar h(z)}=\frac{z}{\sum_{k\geq 1} \frac{z^k}{k!}}=\frac{1}{\sum_{k\geq 0} \frac{z^k}{(k+1)!}}.
\]
Solve the above equation to obtain 

\[
a_0=1, \,\, a_1=-\frac{1}{2}, \,\, a_2=\frac{1}{12},\,\, a_3=0,\,\,  a_4=-\frac{1}{720}\,\, etc., 
\]
which brings us the Riordan array 

\[
[d_{n,k}]_{0\leq k\leq n}=\left[ \frac{k!}{n!} s(n,k)\right]_{0\leq k\leq n}=\left[ \frac{k!}{n!} s(n,k)\right]_{0\leq k\leq n}=\left[ \begin{array}{rrrrr} 1 &  & & &  \\
0 &1 & & & \\ 0 &-\frac{1}{2}& 1& &  \\0 & \frac{1}{3}& -1 &1 & \\ 0 & -\frac{1}{4}& \frac{11}{12}& -\frac{3}{2}& 1\end{array}\right].
\]
The Riordan Stirling array of the signed Stirling numbers of the first kind is 

\[
\left[ s(n,k)\right]_{0\leq k \leq n}= \left[ \frac{k!}{n!} s(n,k)\right]_{0\leq k\leq n}=\left[ \begin{array}{rrrrr} 1 &  & & &  \\
0 &1 & & & \\ 0 &-1& 1& &  \\ 0 & 2& -3 &1 & \\ 0 & -6& 11& -6& 1\end{array}\right],
\]
which gives all $s(n,k)=S(n,k,1,0,0)$ for $0\leq k \leq 4$. For instance, $s(4,0)=0$, $s(4,1)=-6$, $s(4,2)=11$, $s(4,3)=-6$, and $s(4,4)=1$. Of course, the Stirling numbers of the first kind can be evaluated more easily by using formula (\ref{eq:1.15-4}) in Theorem \ref{thm:2.1}, namely,

\[
s(n,k)\equiv S(n,k,1,0,0)=\frac{1}{k!} \left. \frac{d^k}{dz^2} [ z]_{n}\right|_{z=0},
\]
which are simply the coefficients of the powers of $z$ in the expansion of $[z]_n$.

\end{document}